\documentclass[11pt]{article}

\usepackage[applemac]{inputenc}
\usepackage[psamsfonts]{amssymb}
\usepackage[dvips]{graphicx}
\usepackage[psamsfonts]{amssymb}
\usepackage[all]{xy}
\usepackage{lscape}

\def\hH{{H\!H}} % Hochchild H
\def\hC{{\bf C}} % Hochchild C
\def\bB{I\!\!B} % bar B
\def\bH{I\!\!H} % bar H
\def\bk{I\!\!k} % corps
\def\mm{{\mathfrak M}} % model

\title{Rational BV-algebra in String  Topology}
\author{Yves Felix and  Jean-Claude Thomas }

\begin{document}

\maketitle

\rightline{\it To Micheline Vigu\'e-Poirrier on her 60th birthday}
\begin{abstract}
Let $M$ be a 1-connected   closed manifold and $LM$ be the space of
free loops on $M$. In \cite{C-S}  M. Chas and D. Sullivan defined a
structure of BV-algebra on the singular homology of $LM$,  $H_\ast(
LM; \bk)$. When the field of  coefficients is of characteristic zero, we
prove that  there exists a BV-algebra structure on $\hH^\ast(C^\ast
(M); C^\ast (M))$ which carries    the canonical structure of
Gerstenhaber algebra. We construct then an isomorphism  of
BV-algebras between $\hH^\ast (C^\ast (M); C^\ast (M)) $ and the shifted
$ H_{\ast+m} (LM; {\bk} )$. We also prove that the Chas-Sullivan
product and  the BV-operator behave well with the Hodge
decomposition of  $H_\ast (LM) $.
\end{abstract}

\section{Introduction}

 Chas and Sullivan considered in \cite{C-S}  the free loop space $LM= \mbox{map}(S^1, M)$
 for a smooth orientable closed manifold of dimension $m$. They use geometric methods to
  show that the shifted homology $\bH_\ast(LM):=H_{\ast+m}(LM)$ has the structure of a
  Batalin-Vilkovisky algebra (BV-algebra for short).  Later on  Cohen and Jones defined
  in \cite{C-J} a ring spectrum structure  on the Thom spectrum $LM^{-TM}$
  which realizes  the Chas-Sullivan product in homology, that is called the loop product. More recently,
  Gruher and Salvatore proved in \cite{G-S} that the algebra structure
  (and thus the BV-algebra structure) on $\bH_{\ast}(LM)$ is natural with
   respect to smooth orientation preserving homotopy equivalences.

Assume  that the coefficients ring is a field. By a result of
Jones, \cite{J} there exists a linear  isomorphism  $
\hH_\ast(C^\ast(M);C^\ast(M)) \cong  H^\ast(LM)\,,$ and by duality
an isomorphism $ H_\ast(LM) \cong  \hH^\ast(C^\ast(M);C_\ast(M))$.
Here $\hH_\ast (A;Q)$ (respectively  $\hH^\ast (A;Q)$)  denotes
the Hochschild homology (respectively cohomology) of a
differential graded algebra $A$ with coefficients in the
differential graded $A$-bimodule $Q$,   $C^\ast(M)$ denotes the
singular cochains algebra and  $C_\ast(M)$  the complex of
singular chains. The cap product induces  an isomorphism of graded
vector spaces (for instance see \cite{FTV1}-Appendix), $
 \hH ^\ast( C^\ast (M); C_\ast(M))  \cong \hH ^{\ast-m} ( C^\ast (M); C^\ast(M))$,
 and therefore an   isomorphism
 $$
\bH_\ast (LM) \cong \hH^\ast(C^\ast(M);C^\ast(M))\,.$$
Since $ \hH^\ast(A;A)$ is canonically a Gerstenhaber algebra, for any differential graded algebra $A$,
it is natural to ask:

\vspace{3mm} \noindent {\bf  Question 1.} {\it  Does there exist
an  isomorphism of Gerstenhaber algebras   between $\bH_\ast (LM)$
and $ \hH^\ast(C^\ast(M);C^\ast(M))$?}

\vspace{2mm} Various isomorphisms of graded algebras have been
constructed. The first one has been constructed by  Cohen and
 Jones    for any field of coefficients in terms of spectra,
  \cite{C-J}. Then  Merkulov  did it for   for real coefficients \cite{Me},  \cite{Fu}
  using iterated integrals. Finally  for  rational coefficients,  M. Vigu\'e and the two authors
  have constructed another isomorphism, \cite{FTV2}, using the
  chain coalgebra
  of the Quillen minimal model of $M$.

Although $\hH^\ast(A;A)$  does not have, for any differential graded algebra $A
$,  a  natural structure of BV-algebra extending the canonical Gerstenhaber algebra, a
second natural question is:

  \vspace{3mm}
\noindent {\bf  Question 2.} {\it Does there exist on $\hH^\ast(C^\ast(M);C^\ast(M))$ a structure of
BV-algebra extending the   structure of Gerstenhaber algebra  and an isomorphism of BV-algebras between
 $\bH_\ast (LM)$ and $ \hH^\ast(C^\ast(M);C^\ast(M))$?}

\vspace{2mm}
The main result of this paper  furnishes a positive answer to Question 2 and thus to
 Question 1 when the field of coefficients is assumed of characteristic zero.

\vspace{3mm} \noindent {\bf  Theorem 1.} {\it If $M$ is
1-connected and  the field of coefficients has characteristic zero
then there exists a BV-structure on $
\hH^\ast(C^\ast(M);C^\ast(M))$ and an isomophism of  BV-algebras
$\bH_\ast (LM) \cong \hH^\ast(C^\ast(M);C^\ast(M))$.}

\vspace{2mm}   BV-algebra structures on the Hochschild cohomology
$\hH^\ast(A;A)$ have been constructed by different authors under
some conditions on $A$. First of all, Tradler and Zeinalian did it
when $A$ is the dual of an $A_\infty$-coalgebra with
$\infty$-duality (rational coefficients), \cite{T-Z}. This is in
particular the case when $A=C^\ast(M)$, \cite{T}.  In \cite{M1}
Menichi constructed also a BV-structure in the case when $A$ is a
symmetric algebra (any coefficients). Let mention also the
construction due to Vaintrob of a BV-structure on $\hH^\ast(A;A)$
when $A$ is the group ring with rational coefficients of the
fundamental group of an aspherical manifold \cite{Va}.  Vaintrob
gives also an isomorphism of BV-algebras between $\bH(LM)$ and
$\hH^\ast(A;A)$. This is coherent with our Theorem 1  because in
this case $C_\ast(\Omega M) $ is quasi-isomorphic to $A$ and using
\cite{FMT}, we have isomorphisms of Gerstenhaber algebras
$$\hH^\ast(A;A) \cong \hH^\ast(C_\ast(\Omega M);C_\ast( \Omega M))\cong
\hH^\ast(C^\ast(M);C^\ast(M))\,.
$$

Extending  Theorem 1 to finite fields of coefficients would be difficult.
For instance Menichi proved in \cite{M2} that the algebras $\bH(LS^2)$ and
$\hH^\ast(H^\ast(S^2);H^\ast(S^2))$ are isomorphic as Gerstenhaber algebras but
not   as BV-algebras for $\mathbb Z/2$-coefficients.

\vspace{4mm} In this paper we work over a  field of characteristic
zero. We use rational homotopy and we refer to \cite{FHT} for
rational models. We only recall  that   a  morphism  in some
category of complexes is a {\it quasi-isomorphism} if it induces
an isomorphism in homology.  Two objects   are {\it
quasi-isomorphic} if they are related by a finite sequence  of
quasi-isomorphisms.

 Now  by  a recent result of Lambrechts and Stanley, \cite{L-S} there is
 a commutative differential graded algebra $A$ satisfying:
\begin{enumerate}\item  $A$  is quasi-isomorphic to   the differential graded
algebra $C^\ast(M)$. \item  $A$ is connected, finite dimensional
and satisfies Poincar\'e duality in dimension $m$. This means
there exists a $A$-linear isomorphism $\theta : A  \to A^\vee $ of
degree $-m$ which commutes with the differentials.\end{enumerate}
We call $A$ a Poincar\'e duality model for $M$. This model can
replace $C^\ast (M)$ in the Hochschild cohomology because there is
an isomorphism of Gerstenhaber algebra between
$\hH^\ast(C^\ast(M);C^\ast(M))$ and $\hH^\ast(A;A)$ (\cite{FMT}).

Denote by $\mu$  the multiplication of $A$. We define the linear map
 $\mu_A : A \to A\otimes A$    by the commutative diagram
$$
\xymatrix{
\ar@{}[d]^{(\ast)} && A^\vee \ar[r]^(0.3){\mu^\vee}  & (A\otimes A)^\vee =A^\vee  \otimes A^\vee\\
&&A\ar[r]^{\mu_A} \ar[u]^{\theta}  & A \otimes A \ar[u]_{\theta \otimes \theta}
}
$$
 By  definition $\mu_A$ is a degree $m$ homomorphism of $A$-bimodules
 which commutes with the differentials.

 Let  $\hC_\ast (A;A):=(A\otimes T(s\bar A),\partial )$ be the Hochschild
 complex of $A$ with coefficients in $A$.
 Here $T(s\bar A)$ denotes the free coalgebra  generated by
 the graded vector space $s\bar A$ with $\bar A=\{A^i\}_{i\geq 1} $ and $(s\bar A)^i=A^{i+1}$.
 We prove:

 \vspace{2mm}
 \noindent {\bf Proposition 1.}  {\it
 \begin{enumerate}

  \item   The   complex $\hC_\ast (A;A)$ is quasi-isomorphic to $C^\ast(LM)$.

\item  If  $\phi : T(s\bar A)\to T(s\bar A) \otimes T(s\bar A) $ denotes  the coproduct of
$T(s\bar A)$ then the composition
 $$ \Phi :
 A\otimes T(s\bar A) \stackrel{id\otimes \phi }\to  A\otimes T(s\bar A) \otimes T(s\bar A)
 \cong A\otimes_{A^{\otimes 2}} (A\otimes T(s\bar A))^{\otimes 2} \stackrel{\mu_A\otimes id}\to
  A^{\otimes 2}\otimes _{A^{\otimes 2}} (A\otimes T(s\bar A))^{\otimes 2}
 $$
   is a linear map of degree $m$ which commutes with the differentials.

   \item The isomorphism
    $\hH_{\ast}(A;A)^{\ast+m}    \cong   H^\ast(LM) $
     transfers the map induced by $\Phi$ to   the  dual to the  Chas-Sullivan   product.

\item  The duality isomorphism  $ \hH_{\ast}(A;A) ^\vee \cong
\hH^{\ast}(A;A^\vee)  \stackrel{(\theta)}{\cong} \hH^{\ast-m}
(A;A)$
 transfers the map  induced  by $\Phi$ on $\hH_\ast(A;A)$ to the Gerstenhaber
 product on $\hH^\ast(A;A)$.

  \end{enumerate}
 }

 \vspace{3mm}
 Denote by $\Delta : \bH_\ast (LM) \to \bH_{\ast +1}(LM)$ and
 $ \Delta' : \bH^\ast (LM) \to \bH^{\ast -1}(LM) $ the morphisms
  induced by the canonical action of $S^1$ on $LM$.
  As proved by Chas and Sullivan the morphism $\Delta$ together with the loop product equip the
  loop space homology $H_\ast (M)$ with a
   a BV-structure. We  prove:

 \vspace{2mm}\noindent {\bf Proposition 2.}  {\it The  Connes' boundary  $B : \hC_\ast(A;A)
 \to \hC_{\ast -1}(A;A)$    induces  a map in homology corresponding to the operator $\Delta '$  via the
 isomorphism
 $  \hH_\ast (A;A)   \cong   H^\ast(LM) $.}

 \vspace{2mm} On the other hand, the duality isomorphism
 $ \hH_{\ast}(A;A)  ^\vee \cong \hH^{\ast}(A;A^\vee)  \stackrel{(\theta)}{\cong} \hH^\ast (A;A)$
   transfers the linear dual of the Connes'operator $B$  to an operator   on $\hH^\ast(A;A) $
   that gives to the Hochschild cohomology a BV-structure extending the Gerstenhaber
   algebra structure (\cite{M1},\cite{T}) and
  $\hH^\ast(A;A) \cong \hH_\ast(C^\ast (M);C^\ast(M))$ as Gerstenhaber algebras \cite{FMT}.
  This fact combined with Proposition 1 and 2 gives Theorem 1.

   \vspace{3mm}  Since   the field of coefficients is   of characteristic zero,
   the homology of $LM$ admits a Hodge decomposition, $\bH_\ast(LM)   =
   \oplus_{r\geq 0} \bH_\ast^{[r]}(LM)$, \cite{Ge-S},  \cite[Theorem 4.5.10]{Lo}.
   We prove that this decomposition behaves well with respect to the
   loop product $\bullet$   and the BV-operator $\Delta$ defined by Chas-Sullivan.

    \vspace{2mm}\noindent {\bf Theorem 2.}  {\it With the above notations, we have
    \begin{enumerate}
    \item[$\bullet$]  $\bH_\ast^{[r]}(LM)\otimes \bH_\ast^{[s]}(LM) \stackrel{\bullet}{\to}
    \bH_\ast^{[\leq r+s]}(LM)$,
    \item[$\bullet$] $\Delta : \bH_\ast^{[r]}(LM) \to \bH_\ast^{[r+1]}(LM)$\,.
    \end{enumerate}
   }

 \vspace{2mm} By definition  $\bH_\ast^{[0]} (LM)$ is the image    of $H_{\ast -m}(M) $
  by the homomorphism induced in homology by the canonical  section $M\to LM$. It has
  been proved in  \cite{FT}  that  if $\mbox{aut}M$ denotes  the monoid of (unbased)
  self-equivalences of $M $ then there  exist an isomorphism of graded algebras
  $ \bH_\ast^{[1]}(LM) \cong H_{\ast -m}(M) \otimes \pi_\ast (\Omega aut M)$. For any $r \geq 0$,
   a description of
 $\bH_\ast^{[r]}(LM)$ can be obtained,
  using  a  Lie model $(L,d)$ of $M$,  as proved in the last result.

\vspace{2mm}\noindent {\bf Proposition 3.} {\it The graded vector space  $\bH_\ast^{[r]}(LM)$
is isomorphic to $\mbox{\rm Tor}^{UL}(\bk, \Gamma^r(L))$ where $\Gamma^r(L)$ is
 the sub $UL$-module of $UL$ for the adjoint representation that is the image of $\bigwedge^r L$
 by the classical Poincar\'e-Birkoff-Witt isomorphism of coalgebras $\land L \to UL$.}

\vspace{3mm} The text is organized as follows. Notation and
definitions are made precise  in  sections 2 and 3. Proposition 1
is proved in  Sections 4, Proposition 2 is proved in section 5.
Theorem 2  and Proposition 3 are  proved  in the last section.

\section{ Hochschild homology and cohomology}

 \subsection {Bar construction}

  Let  $A$ be a differential graded
augmented
 cochain algebra and let $P$ (respectively   $N $) be a
 differential graded  right (respectively left) $ A $-module, $
 A =   \{A^i\}_{i\geq 0}\,   P =\{P^j\}_{j \in {\mathbb Z}}\,,  N =\{N^j\}_{j \in {\mathbb Z}}$
and $\bar A =$ ker$(\varepsilon : A \to \bk )$. The {\it two-sided
(normalized)  bar construction},
$$
 {\bB} (P; A; N) =  P \otimes T (s  \bar A ) \otimes N\,,
\quad {\bB} _k (P; A; N) =  P \otimes T^k (s  \bar A )
\otimes N\,,
 $$
is the cochain complex defined as follows: For $k\geq 1$, a generic element    $
 p[a_1|a_2|...|a_k]n $ in $
 {\bB} _k (P; A; N)
 $  has  (upper) degree    $| p|+ | n|+ \sum_{i=
 1}^k (|s a_i|)\, $.  If $k=0$, we write $p[\,]n = p\otimes 1\otimes n \in
P \otimes T^0(s\bar A)\otimes N$. The
 differential $d= d_0+d_1$ is defined by:
 $$
  \renewcommand{\arraystretch}{1.6}
  \begin{array}{l}
   {\bB} _k (P; A; N)^{l} \stackrel{d_0}\to  {\bB} _{k} (P; A; N)^{l+1} \,, \quad {\bB} _k
   (P; A; N)^{l} \stackrel{d_1}\to  {\bB} _{k-1} (P; A; N)^{l+1} \\
 \begin{array}{rll}
   d_0  ( p[a_1|a_2|...|a_k]n)& = d( p) [a_1|a_2|...|a_k]n  -
\displaystyle\sum
 _{i=1}^k (-1)^{\epsilon _i}  p[a_1|a_2|...|d(
 a_i)|...|a_k]n\\ & \hspace{5mm} + (-1) ^{\epsilon _{k+1}}  p[a_1|a_2|...|a_k]d(n)
 \\[3mm]
 d_1  (p[a_1|a_2|...|a_k]n) &= (-1) ^{| p|}  pa_1[a_2|...|a_k]n +
 \displaystyle\sum _{i=2}^k (-1) ^{\epsilon _i}
 p[a_1|a_2|...|a_{i-1}a_i|... | a_k]n \\
 &\hspace{5mm} - (-1)^{\epsilon _{k}}  p[a_1|a_2|...|a_{k-1}]
 a_k n
 \end{array}
 \end{array}
 \renewcommand{\arraystretch}{1}
 $$
  Here
 $\epsilon _i = | p| + \sum _{j<i} (|s a_j|)$.

In particular, considering $\bk $ as a trivial $A$-bimodule we obtain  the complex
  $
 \bB A ={\bB} (\bk ;A;\bk )
 $ which  is a differential graded coalgebra whose comultiplication is defined by
 $$\phi ([a_1\vert \cdots \vert a_r] )=
 \sum_{i=0}^r \, [a_1\vert \cdots \vert a_i]\otimes [a_{i+1}\vert \cdots \vert \vert a_r]\,.$$

 Recall that a differential $A$-module $N$ is called { \it semifree}
if $N$ is the union of an increasing sequence of sub-modules
$N(i)$, $i\geq 0$, such that each $N(i)/N(i-1)$ is an $R$-free
module on a basis of cycles (\cite{FHT}). Then,

 \vspace{3mm}
 \noindent {\bf Lemma 1.}
 \cite[Lemma 4.3]{FHT} {\sl    The canonical map  $
 \varphi :  \mathbb B(A;A;A)\to A$  defined by $\varphi [\,] = 1$
 and $\varphi([a_1\vert \cdots \vert a_k]) = 0$ if $k>0$,
 is a
semifree resolution of $A$ as an $A$-bimodule.}

\subsection{Hochschild complexes}

Let us denote by $A^e=A \otimes A^{op}$ the envelopping algebra of $A$.

If  $P$ is a differential graded (right) $A$-bimodule then the cochain complex
$$\hC_\ast(P;A) := (P\otimes T(s\bar A), \partial)  \stackrel{def}\cong
 P \otimes _{A^e} \bB(A;A;A)   $$
  is called the {\it Hochschild chain complex of $A$ with coefficients in   $P$}.
 Its homology is called the {\it Hochschild homology of $A$ with
coefficients in $P$} and is denoted by $
 \hH_\ast(A; P)$. When we consider $\hC_\ast(A;A)$ as well as $\hH_\ast (A;A)$,   $A$ is supposed
  equipped with its canonical   bimodule structure.

\vspace{3mm}
  If $N$ is a (left) differential graded $A$-bimodule then  the cochain complex
  $$\hC ^\ast(A; N) :=  (\mbox{Hom}(T(s\bar A),N), \delta ) \stackrel{def}\cong
  \mbox{Hom}_{A^{e}}(\bB(A;A;A),N)$$
 is called the {\it Hochschild cochain complex} of $A$ with coefficients in the differential
 graded $A$-bimodule $N$.
 Its cohomology is called the {\it Hochschild cohomology of $A$ with
coefficients in $N$} and is denoted by
 $\hH^*(A;N)$.
 When we consider $\hC^\ast(A;A)$ as well as $\hH^\ast (A;A)$,   $A$ is supposed equipped with its
 canonical   bimodule structure.

\vspace{3mm}
Let us denote by $V^\vee $ the {\it graded dual } of the graded vector space
$V =\{V^i\}_{i \in \mathbb Z}$, i.e.
  $V^\vee=\{V^\vee _i\}_{i\in \mathbb Z} \mbox{ with }
V^\vee _i :=    \mbox{Hom}(V^{i}, \bk ) \,.$
The canonical isomorphism
$$
\mbox{Hom}( A \otimes _{A^e} \bB(A;A;A), \bk) \to \mbox{Hom}_{A^e}( \bB(A;A;A), A^\vee)$$
induces the isomorphism  of complexes
 $
\hC_{\ast}(A;A)^\vee \to \hC^\ast(A;A^\vee)$.

\subsection{The   Gerstenhaber algebra on $\hH^\ast (A;A)$}

Ê{\it A Gerstenhaber algebra} is a commutative graded  algebra
$H=\{H_i\}_{i\in \mathbb Z}$  with  a  bracket
 $$
 H_i \otimes H_j \to H_{i+j-1} \,, \quad x\otimes y \mapsto \{x,y\}
 $$
  such that for $a,a', a'' \in H$: \\
(a) $\{a,a'\} = (-1)^{(|a|-1)(|a'|-1)} \{a',a\}$\\
(b) $ \{a,\{a',a''\}\} = \{\{a,a'\},a''\} + (-1)^{(|a|-1)(|a'|-1)} \{a',\{a,a''\}\} $.

For instance the Hochschild cohomology   $\hH^\ast  (A;A) $ is a
Gerstenhaber algebra \cite{Ge}.  The bracket can be defined
   by identifying  $\hC^\ast (A;A)$ with
    a differential graded Lie algebra of coderivations  (\cite{Sta}, \cite[2.4]{FMT}).

\subsection{  BV-algebras and  differential graded Poincar\'e duality algebras.}

A Batalin-Vilkovisky   algebra  (BV-algebra  for short)  is
 a  commutative graded algebra, $H$ together  with a linear map (called a BV-operator)
$$
\Delta : H^k \to H^{k-1}
$$
such that\\
(1) $\Delta\circ \Delta =0$\\
(2) $H$ is a Gerstenhaber algebra with  the bracket  defined by
$$
 \{a,a'\}= (-1)^{|a|}\left( \Delta(aa') -\Delta(a)a' -(-1)^{|a|} ab\Delta(a') \right)\,.
$$

 A differential graded Poincar\'e algebra $(A,d)$ is a finite dimensional commutative differential
 graded algebra together with an isomorphism of differential $A$-modules,
 $\theta : A \to A^\vee$. When $(A,d)$ is a differential graded Poincar\'e algebra,
  then $\hH^\ast (A;A)$ is a BV-algebra \cite{M1}. The BV-operator $\Delta$ is obtained from
   the Connes' boundary  on $\hC_\ast (A;A)$:
  $$B : C_n(A;A)\to C_{n+1}(A;A)\,,$$
 $$
 B(a_0\otimes [a_1|...|a_n]) =
 \left\{
 \begin{array}{ll} 0 & \mbox{ si } |a_0|=0 \\
 \sum _{i=0}^n (-1)^{\bar \epsilon_i}1\otimes  [a_i|...|a_n|a_0|a_1|...|a_{i-1}] &\mbox{ si } |a_0|>0
 \end{array}
 \right.
$$
where
$$
\bar \epsilon_i =(|sa_0|+|sa_1|+...+|sa_{i-1}|)(|sa_i|+...+|sa_n|)\,.
$$
It is well known that $ B^2= 0$ and $ B\circ \partial + \partial \circ B =0\,.$

The operator $\Delta$ is   the image of $B^\vee$ by the duality
isomorphism  $  (\hH_{\ast+m}(A;A))^\vee \cong
\hH^{\ast+m}(A;A^\vee)\stackrel{(\theta)}{\cong} \hH^\ast (A;A)$.

  \section{The Chas-Sullivan   algebra structure on $\bH_\ast(LM)$ and its dual}

    We assume in this section and in the following ones that $\bk $ is a field of characteristic zero.

Denote by $p_0 : LM \to M$ the evaluation map at the base point,
and recall that the space $LM$      can be   replaced by a smooth
manifold (\cite{Ch}, \cite{St}) so that $p_0 $ is a smooth locally
trivial fibre bundle (\cite{Bry}, \cite{St}).

 The loop product
 $$
\bullet : H_\ast( LM  )^{ \otimes 2} \to H_{\ast-m}(LM  ) \,,\quad  x\otimes y
\mapsto
x\bullet y $$
 was first defined by M. Chas and D. Sullivan, \cite{C-S},  by  using  ``transversal geometric
 chains''. With the loop product,
  $\bH_\ast(LM ) := H_{\ast +m} (LM )$ is a
commutative graded algebra.

  It  is  convenient
for our purpose to  introduce the {\it dual of the loop
product}
$ H^\ast( LM ) \to H^{\ast +m}(LM^{\times 2} )
 $. Consider the commutative diagram
$$
(1) \qquad \qquad
\begin{array}{ccccccl}
 LM ^{\times 2} &\stackrel{i }{\longleftarrow}&  LM \times_M LM&
\stackrel{\mbox{\scriptsize Comp}}{\longrightarrow}&  LM \\
\hspace{- 9mm} {\scriptstyle p_0^{\times 2}} \downarrow &&
\hspace{- 4mm}{\scriptstyle p_0} \downarrow && \hspace{2mm}
\downarrow {\scriptstyle p_0}
\\
M^{\times 2}  &\stackrel{\Delta}{\longleftarrow} &M &=& M
\end{array}
$$
where\\
- $\mbox{Comp}$ denotes composition of free loops,\\
- the left hand square is a pullback diagram of locally trivial fibrations,\\
-   $ i$ is the   embedding of the
manifold of composable loops into the product   $LM\times LM$.\\

The embeddings $\Delta$ and $i$ have both codimension $m$. Thus,
using the   Thom-Pontryagin construction
   we obtain the Gysin maps:
 $$
  \Delta ^! : H^k ( M)
\to H^{k+m}(M^{\times 2}) \,, \qquad i^! : H^k (LM\times _MLM) \to
H^{k +m}(LM ^{\times 2}) \,. $$
Thus diagram (1) yields the following diagram:
$$
(2) \qquad \qquad
\begin{array}{ccccccl}
 H^{k+m} (LM ^{\times 2}) &\stackrel{i^! }{\longleftarrow}&  H^{k }(LM
\times_M LM)& \stackrel{H^k(\mbox{\scriptsize
Comp})}{\longleftarrow}&  H^{k } (LM)
\\ \hspace{- 9mm} {\scriptstyle H^\ast(p_0)^{\otimes 2} } \uparrow && \hspace{-
4mm}{\scriptstyle H^\ast(p_0) } \uparrow && \hspace{2mm}\uparrow
{\scriptstyle H^\ast(p_0)}
\\
H^{k+m} (M^{\times 2})   &\stackrel{\Delta ^!}{\longleftarrow}
&H^{k} (M) &=& H^k(M )
\end{array}
$$

\vspace{2mm}\noindent
 Following (\cite{S2},  \cite{C-J-Y}), the  {\it dual of the loop
 product}
is defined by composition of maps on the upper line :
 $$ i^! \circ
H^\ast (\mbox{\scriptsize  Comp})  : H^\ast(LM) \to  H ^ {\ast+m}
(LM  ^{\times 2})\,. $$

 \section{A Hochschild chain model for the dual of the loop product}

  The composition of free loops  $Comp : LM\times _MLM \to LM  $ is obtained by pullback
   from  the composition of paths $Comp' : M^I\times _{M} M^I \to M^I $
   in the following commutative diagram.
  $$ \xymatrix{
\ar@{}[ddd]^{(Comp)} &&&&LM\times_M LM \ar[rr]^{j}
\ar'[d][dd]_{ev_0}\ar[dl]^{\mbox{\scriptsize Comp}} && M^I \times
_M M^I \ar[dd]^{ev_0,ev_1,ev_0} \ar[dl]^{\mbox{\scriptsize Comp'}}
\\
 &&& LM \ar[rr]^(0.6){j } \ar[dd]_{ev_0=p_0} && M^I \ar[dd]^(0.3){(ev_0,ev_1)}&& \\
&&&&M \ar'[r][rr]_(-0.3){(id\times \Delta)\circ \Delta}
\ar@{=}[dl]
&& M^{\times 3}\ar[dl]_{pr_{13}} \\
&&&M\ar[rr]_{\Delta} &&M^{\times 2} }
$$
Here $\Delta$ denotes the diagonal embedding, $j$ the obvious
inclusions,  $ev_t$ denotes the evaluation maps at $t$, and
$pr_{13}$ the map defined by  $pr_{13}(a,b,c)=(a,c)$.

 Let $(A,d)$ be  a  commutative differential graded algebra quasi-isomorphic to the differential graded algebra $C^\ast(M)$.
 A cochain model of the right hand square in diagram $(Comp)$  is
 given by the diagram
$$
\xymatrix{
\ar@{}[d]^{(\dagger)} &
\bB(A;A;A) \ar[r]^(0.3){\Psi }& \bB(A;A;A)\otimes _A\bB(A;A;A)\\
&A^{\otimes 2}\ar[u] \ar[r]^{\psi} &A^{\otimes 3}\ar[u]
}
$$
where $\Psi $ and $\psi$ denote the homorphism  of cochain
complexes defined by
$$
\Psi(a\otimes[a_1|...|a_k]\otimes a')= \sum_{i=0}^k
a\otimes[a_1|...|a_i]\otimes 1 \otimes[a_{i+1}|...|a_k]\otimes a'
\mbox{, and } \psi ( a\otimes a') = a\otimes 1\otimes a'\,.
$$
 The multiplication $\mu : A\otimes A \to A$ makes $A$ a $A$-bimodule, and is a model for the diagonal
  $\Delta$. We consider now the diagram obtained by tensoring diagram $(\dagger)$ by $A$ over $A
^{\otimes 2}$.
$$
\xymatrix{
\ar@{}[d]^{(\ddagger)} &
A \otimes _{A^{\otimes 2}} \bB(A,A,A)  \ar[r]^(0.4){id \otimes \Psi }& A
\otimes _{A^{\otimes 2}} \bB(A;A;A)\otimes _A\bB(A;A;A)\\
& A \otimes _{A^{\otimes 2}}A^{\otimes 2}\ar[u] \ar[r]^{id \otimes \psi} &A \otimes _{A^{\otimes 2}}
 A^{\otimes 3}\ar[u]
}
$$
  Since $\bB(A;A;A)$ is a semifree model of $A$ as
  $A$-bimodule, we deduce from \cite{FHT} that  this diagram is a cochain model of the left hand square.
  Thus we have proved:

\vspace{3mm} \noindent {\bf Lemma 2.} {\it  The cochain complex
$\hC_\ast(A;A)$ is a cochain model of $LM$, and we have an
isomorphism of graded vector spaces
$$
\hH_\ast(A,A)\cong H^{\ast}(LM) \,.$$ Moreover, if   $\phi $
denotes the the coproduct of the   coalgebra   $T(s\bar A)$ then the composition $\Phi$, 
  $$
 \hC_\ast(A;A) \cong A\otimes T(s\bar A) \stackrel{id \otimes \phi} \to A \otimes T(s\bar A)
  \otimes T(s\bar A) =\hC_\ast(A; A)\otimes _A\hC_\ast (A;A)\,,
$$
is   model of the composition of   free loops.}

\vspace{3mm}  Let $(A,d)$ be  a differential graded Poincar\'e duality model of $M$, as defined in the Introduction.

Recall now that the Gysin map $\Delta ^!$  of the
diagonal embedding $\Delta  : M \to M\times M$  is the Poincar\'e
dual of the homomorphism $ H_\ast( \Delta) $. This means that the
following diagram  is commutative
$$
\xymatrix{
H_\ast(M) \ar[r]^{H_\ast(\Delta)}& H_\ast(M\times M) \\
H^{\ast } (M) \ar[u]^{-\cap[M]}_{\cong}  \ar[r]^(0.4){\Delta ^!} &
H^{\ast } (M\times M) \ar[u] _{-\cap[M\times M ]}^{\cong} }
$$
Thus  the linear map of degree  $\mu_A = A \to A\otimes A $
defined in the introduction (Diagram $(\ast)$)  is a cochain model
for $ \Delta^!$.

Remark, \cite{Sta},  that we can choose the pullback of a tubular neighborhood
of the diagonal  embedding $\Delta$ as a tubular neighborhood of
the embedding
 $ i : LM\times _LM \to LM \times LM$. We deduce   that

 \vspace{3mm}
\noindent {\bf Lemma 3.} {\it The linear map  of degree $m$
$$
C_\ast(A;A)\otimes _AC_\ast(A,A) \cong  A \otimes_{A^{\otimes 2}}
C_\ast(A;A)^{\otimes 2} \stackrel{\mu_A\otimes id }\to  C_\ast(A;A)^{\otimes 2}
$$
commutes with the differential and induces $i^!$  in homology.}

\vspace{3mm} The next Lemma follows from standard computation. In
combination with Lemma 2 and 3, it gives the proof of Proposition
1 of the Introduction.

\vspace{2mm} \noindent {\bf Lemma 4.} {\it
 The duality isomorphism  $ (\hH_{\ast+m}(A;A))^\vee \cong \hH^{\ast+m}(A;A^\vee)\stackrel{(\theta)}{\cong} \hH^\ast (A;A)$
 transfers the map  induced  by $\Phi$ on $\hH_\ast(A;A)$ to the Gerstenhaber product on $\hH^\ast(A;A)$.}

\vspace{2mm} \noindent{\bf Remark.}  By putting $F_p := A \otimes
\left(T(s\bar A)\right)^{\leq p}$, for $p\geq 0$, we define a
filtration
$$
A \otimes T (s\bar A)
 \supset ... \supset F_p\supset F_{p-1}\supset ... \supset A=F_0
  $$
such that $\partial F_p \subset F_p$ and $\Phi(F_p) \subset
\displaystyle \oplus_{k+\ell =p} F_k\otimes F_\ell$.  The
resulting spectral sequence
$$
E_2^{p,q}= \ H^q(M) \otimes H^p(\Omega M) \Longrightarrow
 H^{p+q}(LM)
$$is the multiplicative Serre spectral sequence for the fibration
$p_0 : LM \to M$. It
 dualizes into a spectral sequence of algebras
 $$
   H_{q+m}(M) \otimes H_p(\Omega M) \Longrightarrow \bH_{p+q}(LM)\,.
 $$
We recover in this way over $\mathbb Q$ the spectral sequence
defined previously by Cohen, Jones and Yan \cite{C-J-Y}.

\section{The canonical circle action on $LM$}

Let $\rho : S^1 \times LM \to LM  $ be the canonical action of the
circle on the space $LM$. The action $\rho$ induces an operator
$\Delta : \bH_\ast(LM)\to \bH_{\ast+1}(LM)$. The loop product
together with $\Delta$ gives to $\bH_\ast (LM)$ a BV-structure
\cite{C-S}.

Denote by $\mm_M=(\bigwedge V,d)$ a (non necessary minimal)
Sullivan model for $M$ \cite[$\S$-12]{FHT2}.  We put $sV=\bar V$
and denote by $S$  the derivation of $\bigwedge V \otimes
\bigwedge \bar V$ defined by $S(v)=\bar v$ and $S(\bar v)=0$ for
$v\in V$ and $\bar v \in \bar V$. Then a Sullivan model for $LM$
is given by the commutative differential graded algebra
$(\bigwedge V \otimes
 \bigwedge \bar V, \bar d)$ where
  $ \bar d(\bar v) =-S(dv)$  \cite{S-V}.
Moreover in \cite{B-V} Burghelea and Vigu\'e prove that a Sullivan
model of the action $\rho: S^1\times LM\to LM$ is given by
$$
\begin{array}{ccc}
\mm_\rho :   (\bigwedge V\otimes \bigwedge \bar V,\bar d) \to
 (\bigwedge u,0) \otimes (\bigwedge V\otimes \bigwedge \bar V,   \bar d)\,,\hspace{5mm} | u|=1 ) \\
\mm_\rho (\alpha  ) = 1 \otimes \alpha  +  u\otimes  S(\alpha)\,,
\quad \alpha \in \bigwedge V\otimes \bigwedge \bar V
\end{array}
$$

In particular the map induced in cohomology by the action of $S^1$ on $LM$ is given by the
derivation $S : H^\ast (\bigwedge V\otimes\bigwedge \bar V)\to
 H^{\ast -1}(\bigwedge V\otimes\bigwedge \bar V)$. Denote now by $B$ the Connes'
 boundary on $\hC^\ast (\mm_M;\mm_M)) =\bigwedge V \otimes T(s\overline{\bigwedge V})$.
 Vigu\'e proved the  following Lemma  in  \cite[Theorem 2.4]{Vi}.

\vspace{3mm}\noindent {\bf Lemma 5.}  {\sl The morphism $f :
\hC^\ast (\mm_M;\mm_M) \to
 (\mm_M\otimes\bigwedge \bar V)$ defined by
$$f(a\otimes [a_1|\cdots |a_n])= \frac{1}{n!} a S(a_1)\cdots S(a_n)$$
is a quasi-isomorphism of complexes and $f\circ B = S \circ f$.}

\vspace{3mm} This, combined with lemma 4, implies directly
Proposition 2 of the Introduction.

\section{The Hodge decomposition}

With the notation of the previous sections, let  $(\bigwedge V
\otimes \bigwedge \bar V, \bar d)$ be a Sullivan model for $LM$.
Denote by $G^p = \bigwedge V \otimes \bigwedge ^p \bar V$ the
subvector space  generated by the words of length $p$ in $\bar V$.
The differential $ \bar d$ satisfies
 $ \bar d ( G^p  ) \subset  G^p$.
Thus we put
$$
H^n_{[p]} ( {LM}) := H^n ( G^p) \,.
$$
This decomposition has been considered by many authors (see for
instance \cite{Lo},
 \cite{V}). It induces by duality a Hodge decomposition on
 $H_\ast(LM)$. We are now ready to prove Theorem 2 of the
 Introduction.

\vspace{2mm}\noindent{\bf Proof of Theorem 2.}  Recall that the
differential $d$ in $\hC^\ast (\mm_M;\mm_M)$ decomposes into $d =
d_0 + d_1$ with $d_0( \mm_M\otimes T^p (s\overline{\land V}))
\subset  \mm_M\otimes T^p (s\overline{\land V})$, and $d_1(
\mm_M\otimes T^p (s\overline{\land V}) )\subset  \mm_M\otimes
T^{p-1} (s\overline{\land V})$.

We  consider the quasi-isomorphism  $f : \hC^\ast (\mm_M;\mm_M) \to
(\mm_M\otimes\bigwedge \bar V, \bar d)$ defined in Lemma 5. If we
apply Lemma 5, when $d=0$ in $\land V$, we deduce that $\mbox{Ker}\,
f$ is $d_1$-acyclic. Denote  $K^{(p)}= \mbox{Ker}\, f\cap \left(
\mm_M\otimes T^p (s\overline{\land V})\right)$. Then

\vspace{3mm}\noindent {\bf Lemma 6.}  {\sl  \begin{enumerate}\item
If $\omega \in K^{(p)} \cap \mbox{\rm Ker}\, \partial$ then there
exists $\omega' \in \oplus_{r\geq p+1}K^{(r)}$ such that $\partial
\omega'=\omega $. \item    $f$ induces a surjective map
 $
 \left(\mm_M\otimes T^{\geq p} (s\overline{\land V})\right)  \cap \mbox{Ker}\, \partial
 \twoheadrightarrow  \left(\mm_M\otimes \bigwedge ^p s V\right) \cap \mbox{Ker}\,\bar d
 \,.$
 \end{enumerate}}

 \vspace{3mm}\noindent
 {\bf Proof.}
  If  $\omega \in K^{(p)} \cap \mbox{\rm Ker}\, \partial$  then $\omega = \partial (u+v)$
  with $u \in K^{(p)}$ and $v \in K^{(\geq p+1)}$. Since $\partial _1u=0$ we have
   $u = \partial \beta_1 $ some $\beta \in K^{(p+1)}$ and thus $\omega-d\beta_1 \in K^{(\geq p+1)}$.
   An induction on $n \geq 1$  we prove that there exists  $\beta_n \in K^{(p+n)}$  such that
     $\omega-d\beta_n \in K^{(p+n)}$. Since $\bigwedge V$ is 1-connected
      $\left( \mm_M\otimes T^{p+n} (s\overline{\land V})\right) ^{|\omega|}=0$ for some integer $n_0$.
       We put $\omega'  =\beta_{n_0}$.

   In order to prove  the second statement, we consider a $\bar d$-cocycle
   $\alpha \in \mm_M\otimes\bigwedge^p  s V $  and we write
   $\alpha = f(\omega)$ for some $\omega \in \mm_M\otimes T^{p} (s\overline{\land V}) $.
   It  follows from the definition of $f$ that
   $\partial \omega  \in K^{(p-1)}$. Thus, by the first statement,  $\partial \omega =\partial \omega'$
   some $\omega' \in K^{(\geq p)}$. Then $\varpi = \omega-\omega'$ is $\partial$-cocycle of $K^{\geq p}$
   such that $f(\varpi )=\alpha$.

   \hfill{$\square$}

To end the proof, let $\alpha \in H^\ast_{[n]}(LM)$, then by Lemma
6, $\alpha$ is the class of  $f(\beta)$ where $\beta \in 
\mm_M \otimes T^{\geq n} (s\overline{\bigwedge V})$. Therefore
$\Phi(\beta)\in    \oplus _{i+j\geq n} \left( \mm_M\otimes T^{  i}
(s\overline{\land V})\right) \otimes \left( \mm_M\otimes T^{  j}
(s\overline{\land V})\right) $  (see Lemma 2).   Now since $f( \mm_M\otimes T^p
(s\overline{\land V})) \subset \mm_M\otimes \bigwedge ^p s V$,
$$[\Phi(\alpha)] \in \oplus_{i+j\geq n} H^\ast_{[i]}(LM)\otimes
H^\ast_{[j]}(LM)\,.$$
 \hfill $\square$

\vspace{3mm} Now, as announced in the Introduction (Proposition 3)
there is an other interpretation of $H^n_{[p]} ( {LM}) $ in terms
of the cohomology of a differential graded Lie algebra.

Let $L$ be a differential graded algebra $L$ such that the cochain
algebra ${\cal C}^\ast(L)$ is a Sullivan model of $M$,
\cite[p.322]{FHT2}. In particular, the homology of the enveloping
universal algebra of $L$, denoted $UL$ is a Hopf algebra isomorphic
to $H_\ast(\Omega M)$. We consider the cochain complex ${\cal
C}^\ast (L; UL^\vee _a) $ of $L$ with coefficients in $UL^\vee$
considered as an $L$-module for the adjoint representation.  We have
shown  (\cite[Lemma 4]{FTV2})   that the natural inclusion ${\cal
C}^\ast(L) \hookrightarrow {\cal C}^\ast (L; UL^\vee _a) $ is a
relative Sullivan model of the fibration $p_0: LM \to M$.  Write
${\cal C}^\ast(L)=(\bigwedge V, d)$, then $V=(sL)^\vee $ and $\bar V
=L^\vee$. There is also (Poincar\'e-Birkoff-Witt Theorem) an
isomorphism of graded coalgebras, \cite[Proposition 21.2]{FHT2}:
$$
\gamma : \bigwedge L \to UL\,, \quad x_1\wedge...\wedge x_k \mapsto \sum_{\sigma \in {\mathfrak S}_k} \epsilon_\sigma \, x_{\sigma(1)}...x_{\sigma(k)}\,.
$$
If we put $\Gamma^p = \gamma ( \bigwedge^pV)$ we obtain the
following isomorphisms of cochain complexes
$$
(\bigwedge V \otimes \bigwedge  \bar V,\bar d)\cong  {\cal
C}^\ast(L;UL^\vee_a)\,,\quad  G^p\cong  {\cal
C}^\ast(L;(\Gamma^p)^\vee)
$$
which in turn induce the isomorphisms:
$$
\bH^\ast (LM)\cong \mbox{Ext}_{UL}(\bk, UL_a^\vee )\,, \quad
\bH_{[p]}^\ast (LM)\cong \mbox{Ext}_{UL}(\bk, \Gamma^p(L)^\vee)\,.
$$
and  by duality,
$$
\bH_\ast (LM)\cong \mbox{Tor}^{UL}(\bk, UL_a  )\,, \quad
\bH^{[p]}_\ast (LM)\cong \mbox{Tor}^{UL}(\bk, \Gamma^p)\,.
$$

\footnotesize{

}

\vspace{1cm}

\hspace{-1cm}\begin{minipage}{19cm}
\small
\begin{tabular}{lll}
felix@math.ucl.ac.be
&jean-claude.thomas@univ-angers.fr\\
D\'epartement de math\'ematique    &
D\'epartement de math\'ematique \\

 Universit\'e Catholique de Louvain   &  Facult\'e des
Sciences \\

 2, Chemin du Cyclotron            & 2, Boulevard
Lavoisier     \\

 1348 Louvain-La-Neuve, Belgium       &
 49045
Angers, France

\end{tabular}

\end{minipage}

\end{document}